\documentclass[a4paper]{article}
\usepackage[utf8]{inputenc}
\usepackage[T1]{fontenc}
\usepackage{tikz}
\usepackage{bm}
\usepackage{amsmath}
\usepackage{amsfonts}
\usepackage{xcolor}
\usepackage{authblk}
\usepackage[version=4]{mhchem}

\usepackage[margin=40pt]{geometry}

\tikzset{every picture/.style={line width=0.75pt}} 

\title{Vanquishing the computational cost of passive gamma emission tomography simulations: a physics-aware reduced order modeling approach}

\author[1]{Nicola Cavallini}
\author[1]{Riccardo Ferretti}
\author[1]{Gunnar Bostrom}
\author[2]{Stephen Croft}
\author[1]{Aurora Fassi}
\author[1]{Giovanni Mercurio}
\author[1]{Stefan Nonneman}
\author[1,3]{Andrea Favalli}
\affil[1]{European Commission, Joint Research Centre, Via Enrico Fermi, Ispra, I -211027 (Va), Italy}
\affil[2]{Lancaster University, Bailrigg, Lancaster, UK}
\affil[3]{Los Alamos National Laboratory, P.O. Box 1663, Los Alamos, NM 87545, USA}

\affil[*]{andrea.favalli@ec.europa.eu}

\newcommand{\fbp}{\mathcal{F}}

\begin{document}

\begin{abstract}

Passive Gamma Emission Tomography (PGET) has been developed by the International Atomic Energy Agency as a way to directly image the spatial distribution of individual fuel pins in a spent nuclear fuel assembly and so determine potential diversion.
 
Constructing the analysis and interpretation of PGET measurements rely on the availability of comprehensive datasets. Experimental data are expensive, limited, and so are augmented by Monte Carlo simulations. The main issue concerning Monte Carlo simulations is the high computational cost to simulate the 360 angular views of the tomography. Similar challenges pervade numerical science.
 
To address this challenge, we have developed a physics-aware reduced order modeling approach. It provides a framework to combine a small subset of the 360 angular views with a computationally inexpensive proxy solution, that brings the essence of the physics, to obtain a real-time high-fidelity solution at all angular views, but at a fraction of the computational cost.
\end{abstract}

\flushbottom
\maketitle
%
%
\thispagestyle{empty}

\section*{Introduction}

Passive Gamma Emission Tomography (PGET) is a measurement technique developed for spent nuclear fuel verification to help meet nuclear non-proliferation and safeguards requirements. It is part of the armory of technical measures put in place by the International Atomic Agency (IAEA), in the framework of the Treaty on the Non-Proliferation of Nuclear Weapons (NPT), to ensure that each Member State complies with the requirements of spent nuclear fuel stewardship and accountancy. The PGET system is designed to perform partial defect and bias defect verification of a spent nuclear fuel assembly while it is shielded underwater prior to being moved to long term storage or into a geological repository. The IAEA approved PGET in 2017 for inspections \cite{8533017, iaea-next-gen-2022, Virta2022ER, PNNL-SA-135401, Chernikova2018, Fang2021}.

The PGET detects the passive gamma-ray emissions from the assembly arising from the build-up of fission products generated during the use in the reactor. The most important gamma-ray signatures come from the fission products 
\ce{^{137}Cs}, \ce{^{134}Cs}, \ce{^{154}Eu}, \ce{^{106}Ru}, and \ce{^{144}Ce}  \cite{FAVALLI2016102}. They are important because they confer information about the initial enrichment, burnup and cooling time of the assembly while being practical to measure. Only a few fission products are produced in sufficient quantity, and have intense yet long-lived and penetrating gamma-ray emissions to be of general interest for nondestructive assay applications. 

The PGET design consists of two highly collimated linear arrays comprising up to 91 cadmium-zinc-telluride (CdZnTe) gamma detectors. The detector pitch is 4 mm. Each detector is behind a tungsten collimator 15 mm wide, 100 mm thick, and 5 mm high in front of the detector, increasing to 70 mm at the exit. The two detector heads are on a rotating disk on opposite sides with a 2 mm offset. This 2 mm stagger is devoted to obtaining 182 data collections, referred to as counts, per each view. The technical details of the PGET detector system are reported in \cite{Belanger-Champagne2019,Backholm2020,Virta2022ER}.

The two detectors arrays rotate continuously around the spent fuel assembly. The number of events are integrated both over energy and rotation windows. Typical values for the energy windows are: <400 keV, 400-600 keV, 600-700 keV, and 700-1200 keV, while rotations are integrated over each of the 360 degrees. The end product of a PGET acquisition sequence is a sinogram for each energy window, a matrix characterized by a number of rows equal to the number of detectors and a number of columns equal to the number of angular views, 360 in a standard acquisition. An accurate sinogram is needed as the input to tomographic reconstruction algorithms aimed at reconstructing spent fuel axial cross-sections in order to detect cases of anomalies or diversion, even to the single rod level (so-called bias defect). 

It is worth mentioning that compared with medical and industrial tomography systems, PGET involves peculiar challenges due to the wide range of gamma activities that characterize the spent nuclear fuel and the high gamma attenuation of the nuclear materials. These parameters are outside the control of the experimenter and are challenging in the case of detecting anomalies in the center of the assemblies \cite{Virta2020}.

Development and evaluation studies of algorithms for the analysis and interpretation of PGET measurements, for example, based on artificial intelligence \cite{Fang2021, lahiri:22:svc-arxiv}, rely on the availability of extensive and comprehensive datasets of experimental and/or simulated data. Experimental data are extremely limited in availability and breadth of potential diversion scenarios, and are also expensive to generate. Realistic simulations of spent nuclear fuel measurements are therefore essential in evaluating the performance of the algorithms behind the interpretation of PGET measurements, especially in robustly identifying potential diversion situations.

Monte Carlo simulations, for example using the Monte Carlo N-Particle
(MCNP) code\cite{mcnp-user-manual}, allow the construction of a realistic and extended dataset of simulated cases and related sinograms; however, the computational cost is high. High-fidelity simulations, in fact, require the treatment of the full physics radiation transport in a detailed geometry of the PGET system. Miller et al. and Wittman et al.  pointed out that a high-fidelity MCNP simulation of the PGET measurement of a fuel assembly, the basic 360 angular views of a sinogram, take about 6-7 days in a cluster composed of 128 nodes and 8192 cores \cite{PNNL-SA-135401, PNNL-SA-13389,JNT-1955-Phase}.

To address this challenge, we have developed a physics-aware reduced order modeling approach. In this approach we combine (1) a limited number of the 360 angular views (limited views tomography) and (2) a real time analytically-built approximated sinogram (based on the Lambert-Beer law) at all the 360 views, that brings the essence of the physics. The computational cost is cut by only simulating a sparse subset of angular views. The computational benefit is proportional to the number of skipped views (e.g., simulating 60 views instead of 360 reduces the cost to about one-sixth). The Physics-Aware Reduced Order Model approach then enables the reconstruction of all the angular views. This is achieved using Proper Orthogonal Decomposition (POD)\cite{TezzeleDemoRozza2019MARINE, OrtaliDemoRozza2020MINE, DemoTezzeleRozza2019CRM, certified_rozza_2015}, a subset of the techniques available in Reduced Order Modeling numerical methods, which identifies the dominant patterns in the data by applying Singular Value Decomposition (SVD) to the matrix of the limited views data. 
We named our method Physics-Aware Proper Orthogonal Decomposition (PA-POD).

The method’s performance is tested and measured against the data released by IAEA under the PGET Tomographic and Analysis Challenge (2019
\cite{iaeauniteweb, iaeafullreport}. 
We used the IAEA field dataset as surrogate of realistic simulations.
The dataset is used as the ground truth to which apply our PA-POD method. We sample limited set of angular views, and we reconstruct the sinogram at each angular view via PA-POD. The goodness of the method is obtained by comparison between ground truth and the reconstruction. Among the collection of IAEA data, three cases are of particular interest because they refer to field scenarios. They are named competition three, four, and five, and related to a VVER (water-water energetic reactor), a PWR (Pressurized Water Reactor), and a BWR (Boiling Water Reactor) fuel assembly, respectively. In the Results section we report the results of the PA-POD method applied to the IAEA PWR assembly case, as PWR nuclear plants are the large majority of nuclear plants in the world. All the results presented refer to the 600-700 keV window: it contains the \ce{^{137}Cs} gamma emission line of 661.7 keV, a primary gamma emission of spent nuclear fuel with a relatively long half-life of 30 years\cite{FAVALLI2016102}.

\section*{Our Approach}
We consider the sinogram matrix $\mathbf{S}$, our ground truth,  as a discrete representation of an unknown continuous solution $\mathbf{s}(y,\theta)$, where $y$ represents the spatial coordinate of the detector and $\theta$ is the angular coordinate.  
POD approximates the solution as 
a linear combination of dominant patterns, or modes, 
${\mathbf{u}_i(y)}$ estimated 
using only a reduced set of angular views. The approximated solution is: 
\begin{equation}
\mathbf{\tilde{s}}(y,\theta)=\sum_{i=0}^{k-1} 
\mathbf{u}_i(y)\, \mathbf{c}_i (\theta),
\label{eq:pod-approx}
\end{equation} 
with $k$ the number of selected modes.
All the possible linear 
combination of the modes $\{\mathbf{u}_i(y)\}$, with $i=0,\ldots,k-1$, construct the so-called POD solution space (or POD space), 
where $k$ is its dimesionality.
In the following steps, we fix the notation and summarize our PA-POD 
approach.

\begin{enumerate}
\item \textbf{Database Creation:}
Take the set of $N_s$ simulated/measured limited angular views 
at the angular coordinates $\theta_n$ and construct the set of 
pairs $\{\mathbf{s}(y,\theta_n),\theta_n\}$, with $n = 0,\ldots,N_s-1$. They are stored 
in a database matrix $\mathbf{\hat{S}}_{N\times N_s}$, where 
$N$ is the number of detectors.

\item \textbf{POD Solution Space Construction:}
Apply Singular Value Decomposition to the database matrix: $\mathbf{\hat{S}} = \mathbf{U}\,\bm{\Sigma}\, \mathbf{V}^*$. Meaning decompose 
the database matrix into two unitary matrices $\mathbf{U}$ and $\mathbf{V}$ and a diagonal matrix $\bm{\Sigma}$. Here $^*$ denotes the complex conjugate transpose. The diagonal elements $\delta_j$ of $\bm{\Sigma}$, or singular values, are non-negative and ordered from the largest to the smallest. 
They are usually normalized in the form $\sigma_j = \frac{\delta_j}{\sum_{l=0}^{N-
1}\delta_l}$. We define information variance 
relative to the $j$-th singular value as the sum 
of all the first  
$\sigma_j$. It is a measure of how much of 
the total data variability is captured by the first $j$ modes \cite{Brunton2022, certified_rozza_2015}. 
The $i$-th column of $\mathbf{U}$ stores the mode 
$\mathbf{u}_i(y)$. 
To construct a suitable POD solution space we 
balance between the number of modes and the information 
variance. This reduces the dimensionality
of the problem and avoids selecting higher modes 
\cite{Brunton2022} which are generally 
associated with noise rather than information.

\item \textbf{Coefficients Estimation:}
We use $\mathbf{U^*}_{k\times N}$, the conjugate transpose of the first $k$ columns of the modes matrix $\mathbf{U}$, to project into the POD space $\mathbf{R}_{N\times 360}$, a computationally inexpensive approximation of the full 
sinogram $\mathbf{S}$. We refer to this
approximation as the Real-Time Approximate Forward Model, shortened as Real-Time Model.
In the PGET case it consists in a Lambert exponential attenuation model applied to a voxelization of the fuel assembly geometry.
It simplifies the gamma photon transport and enables a real time approximation at any possible PGET angular view.
This model was initially proposed by  \cite{Backholm2020}, details on our own implementation are described in section Methods.
As a result we obtain the estimation of the set of  
coefficients at the 360 
views, $\{\mathbf{c}_i(\theta_m)\}$ in equation (\ref{eq:pod-approx}), 
with $i = 0,\ldots,k-1$, and $m = 0,\ldots,359$. In matrix form:  $\mathbf{\tilde{C}}_{k\times 360} = \mathbf{U^*}_{k,N}\, \mathbf{R}_{N\times 360}$. 
The coefficients $\tilde{\mathbf{{C}}}$ are row wise scaled to match $\mathbf{\hat{S}}$
values.

\item \textbf{Solution Evaluation:} The final approximated sinogram is 
represented by equation(\ref{eq:pod-approx}), that can be written in matrix form: $\mathbf{\tilde{S}}=\mathbf{U}\, \mathbf{\tilde{C}}$. 
\end{enumerate}

We define and report the error on the reconstructed image rather 
than on the 
sinogram because in the PGET analysis 
the sinogram quality is judged by the
quality of its reconstruction. As a metric we choose to quantify the error matrix $\mathbf{e}$ of the filtered back projection (FBP, mathematically represented with the operator $\mathcal{F}(\cdot)$) 
applied to the ground truth, to the Real-Time Approximated Forward Model, and to our PA-POD estimation. 
Specifically, we used the ramp-filtered version of FBP \cite{kak2001principles,van2014scikit}. In explicit terms, we have the following pixel wise error map definitions:

\[
\mathbf{e}_{\mathrm{Real-Time\ Model}} = \left|
\frac{\fbp(\mathbf{R})-\fbp(\mathbf{S})}
{\fbp(\mathbf{S})}\right|;\quad
\mathbf{e}_{\mathrm{PA-POD}} = 
\left|\frac{\fbp(\mathbf{\tilde{S}})-\fbp(\mathbf{S})}
{\fbp(\mathbf{S})}\right|.
\]

We apply a mask (matrix) to select the $\fbp(\mathbf{S})$
values that are greater than 15\% of the maximum.
The mask covers the assembly cross section and its premises,
this area is addressed as $\Omega_{\mathrm{tot}}$.


Furthermore, we use the cumulative error distribution as an integral error
measure of maps $\mathbf{e}_{\mathrm{Real-Time\ Model}}$ and 
$\mathbf{e}_{\mathrm{PA-POD}}$. Given an error threshold $\epsilon_{\mathrm{th}}$, we count the 
number of pixels with a smaller or equal error, we label it with $\Omega|_{\epsilon_{\mathrm{th}}}$. 
We define the pixel fraction as the ratio between this area and $\Omega_{\mathrm{tot}}$:
\begin{equation}
\mathrm{pixel\ fraction}=\frac{\Omega|_{\epsilon_{\mathrm{th}}}}{\Omega_{\mathrm{tot}}},
\label{ed:pixel-fraction}
\end{equation}
the greater its value, the better the performance.
The metric we 
use in this work does not require any arbitrary 
filtering as the commonly used Structural Similarity Index
\cite{1284395}.

\section*{PWR Case Results}

\begin{figure}[!ht]
\centering
\input{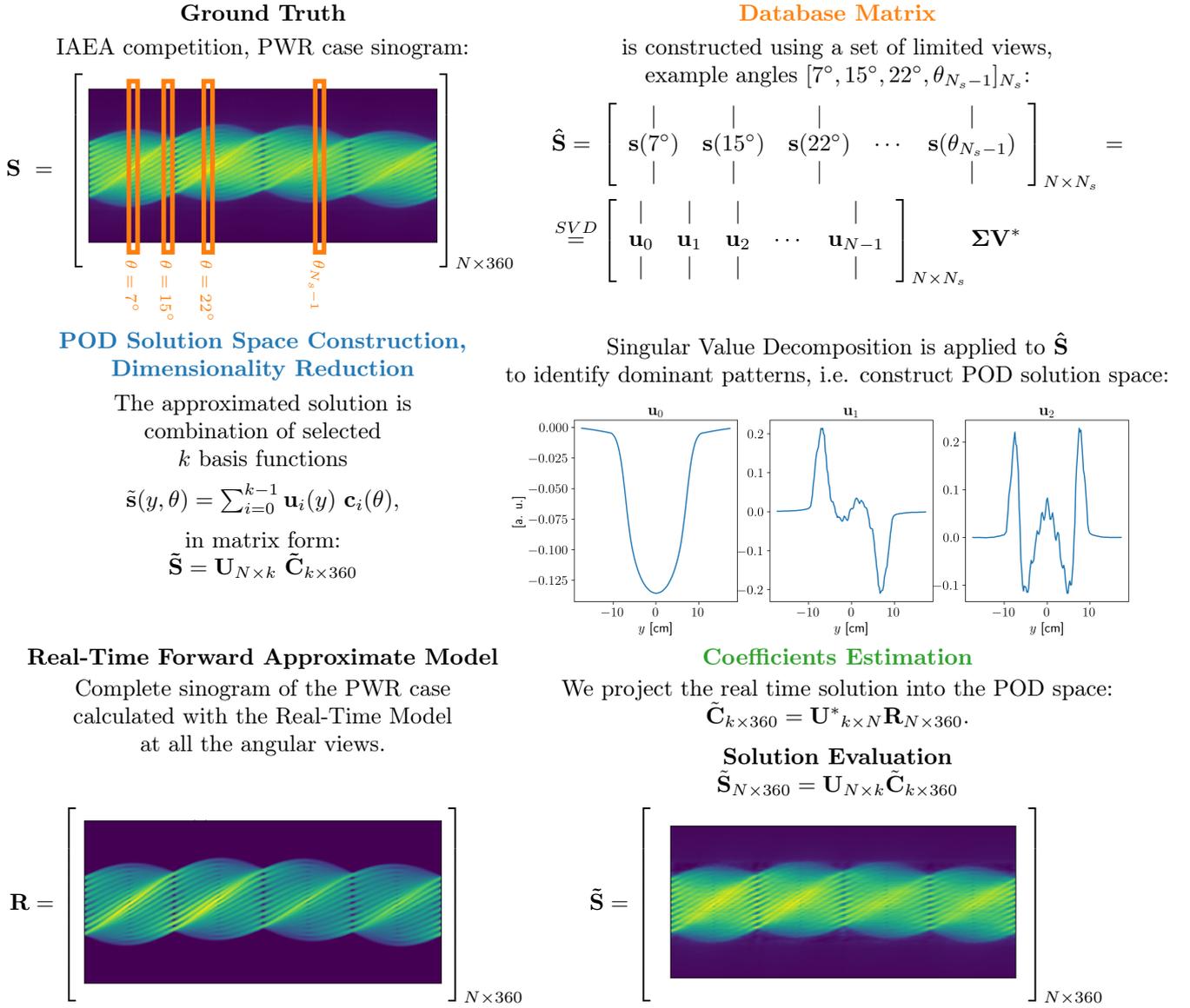}
\caption{The figure reports  the major steps in our PA-POD method. 
On top left we represent 
the actual sinogram for the PWR assembly from the IAEA competition 
(our ground truth), as high fidelity 
data, named $\mathbf{S}$ in our formulation. 
With the orange color we highlight the subset of limited 
views that construct the database matrix. 
We present the POD solution space and plot the first three modes 
$\mathbf{u}_0, \mathbf{u}_1, \mathbf{u}_2$ 
for the case of 
$N_s = k = 60$ randomly chosen views.
On bottom left, the sinogram $\mathbf{R}$, on the same pin 
configuration, as obtained by the Real-Time Approximate Forward Model. 
At the bottom right the PA-POD 
formulation and resulting sinogram $\tilde{\mathbf{S}}$.}
\label{fig:pod-pget-chart}
\end{figure}

The following results report the application of our PA-POD approach, 
starting from a sparse sampling of the 360 angular views of the PWR 
assembly sinogram included in the IAEA competition dataset. 
Specifically, it consists of a 10$\times$10 pins assembly where a 3$\times$3 
block of pins has been removed;
thus it is also a case of a diversion scenario.
The Real-Time Model relies on the investigation domain discretisation, 
see section Methods for the details, 
the reported results are computed with the following 
mesh sizes: $\Delta z =10$ mm, $\Delta x =\Delta y =0.5$ mm,
being $x,y,z$ the three coordinate axes of the 
investigation domain. 
Both the ground truth sinogram 
and the Real-Time one 
are normalized between zero and one.

The flowchart in figure \ref{fig:pod-pget-chart} summarises our method and presents 
an overview of the results of the PWR case. 
In particular, 
we report the ground truth sinogram $\mathbf{S}$, the real-time 
approximated sinogram $\mathbf{R}$, the first three modes $\mathbf{u}_0, \mathbf{u}_1, 
\mathbf{u}_2$ and the PA-POD approximate sinogram $\tilde{\mathbf{S}}$. 
In the real-time approximated forward model, 
we assume the intensity of pins of the 
PWR assembly is flat where the pins are located, 
and zero in the 3 
by 3 block in the lower left center where the pins are removed. 
The PA-POD sinogram $\tilde{\mathbf{S}}$ is constructed randomly 
sampling 60 views 
for the database matrix, and the first 60 modes are selected to 
construct the POD space. The orange views highlighted
in the top left corner have 
the objective of proving a visualisation of the database 
matrix, the rotation angles are 
just examples.

\begin{figure}[h!]
\centering
\includegraphics[width=\linewidth]{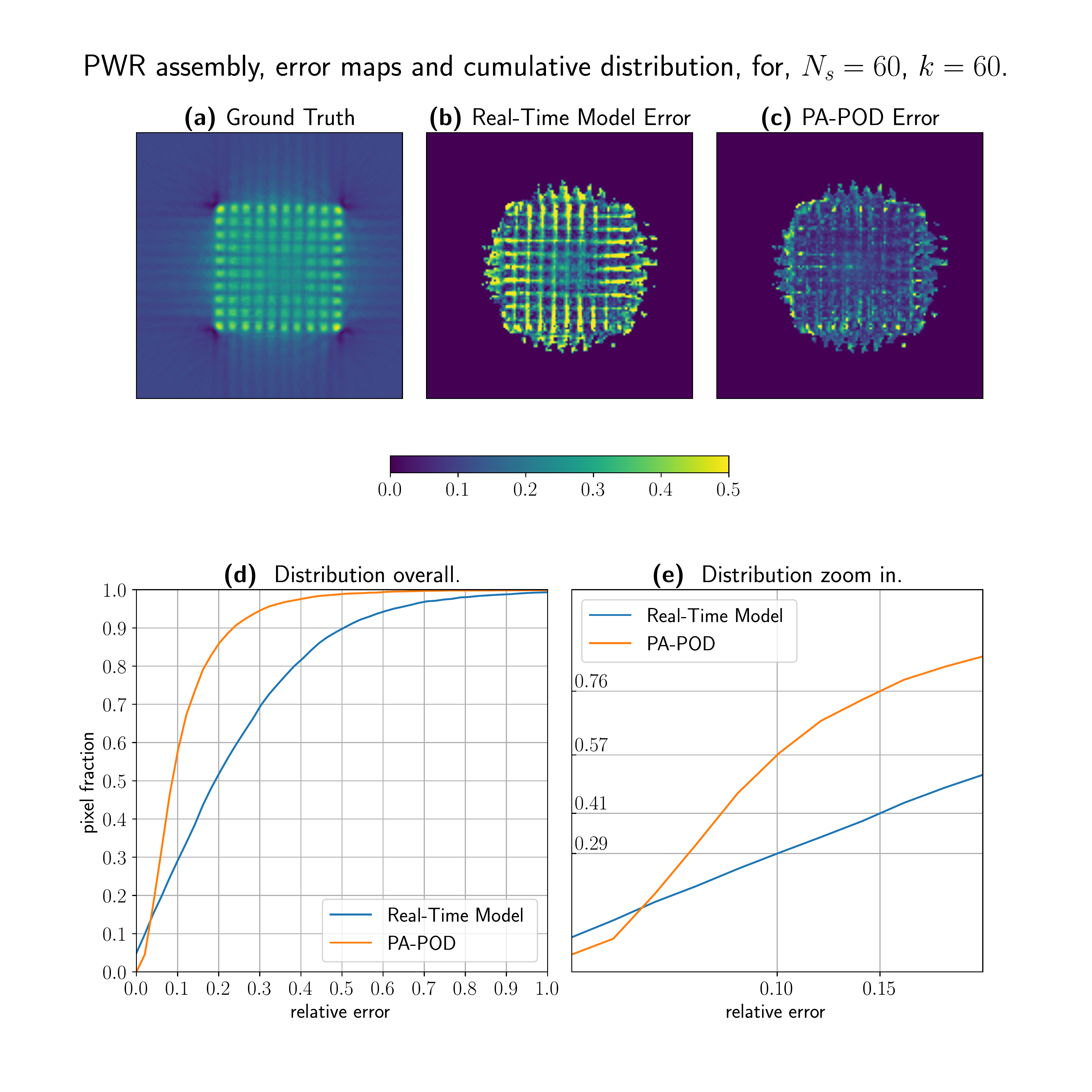}
\caption{Here we represent the ground truth \textbf{(a)} together with the error that characterizes the real time model \textbf{(b)} and our PA-POD approximation \textbf{(c)}. In the PA-POD case we randomly sample sixty views of the spent fuel, we repeat the sample one hundred times, and for each pixel, we collect the median of the sampled data. \textbf{(d)} The improvement from the Real Time Model to PA-POD is quantified by the area between the blue and orange lines. \textbf{(d)} More precisely, considering a 10\% error, PA-POD describes 57\% of the total area, while the Real Time Model describes 29\% of the area.}
\label{fig:c04-pod-relative-error-report}
\end{figure}

Figure \ref{fig:c04-pod-relative-error-report}
represents the error and quantifies it. 
In picture 
\textbf{(a)} 
the ground truth, while the error is spatially 
visualized in the maps \textbf{(b)} and \textbf{(c)}. 
To avoid any procedural 
biasing effect, we repeat the views sampling 100 times, 
and at each pixel we 
plot the median value for the error. 
Figure \textbf{(c)} shows that PA-POD can 
correct most of the high error pixels in between 
pins and around 
the assembly. The PA-POD 
relative error map clearly indicates 
that the approach reproduces the 
actual data within a relative error of the order of 10\%.

Figure \ref{fig:c04-pod-relative-error-report}\textbf{(d)} 
and \textbf{(e)} complete the 
picture plotting the error cumulative distribution, 
that we called pixel fraction to emphasize its geometrical meaning, 
see equation \ref{ed:pixel-fraction}. 
Figure \ref{fig:c04-pod-relative-error-report}\textbf{(d)} starts with 
a minimal number of pixels having the same value as the ground truth due 
to the normalization of the two sinograms. As soon as we consider larger 
errors, the area that PA-POD can describe with a specific accuracy 
increases with a far steeper gradient than the Real-Time Model. Figure 
\ref{fig:c04-pod-relative-error-report} \textbf{(b)}, shows the 
results with a finer-graded detail. We focus on the 10\% relative error 
and observe that PA-POD covers 58\% of the total area, while the Real-
Time Model describes 29\% with the same accuracy. These results show 
that it is possible to  approximately double the accuracy of the real-
time model using 16\% of the total number of views.

\begin{figure}[h!]
\centering
\includegraphics[width=\linewidth]{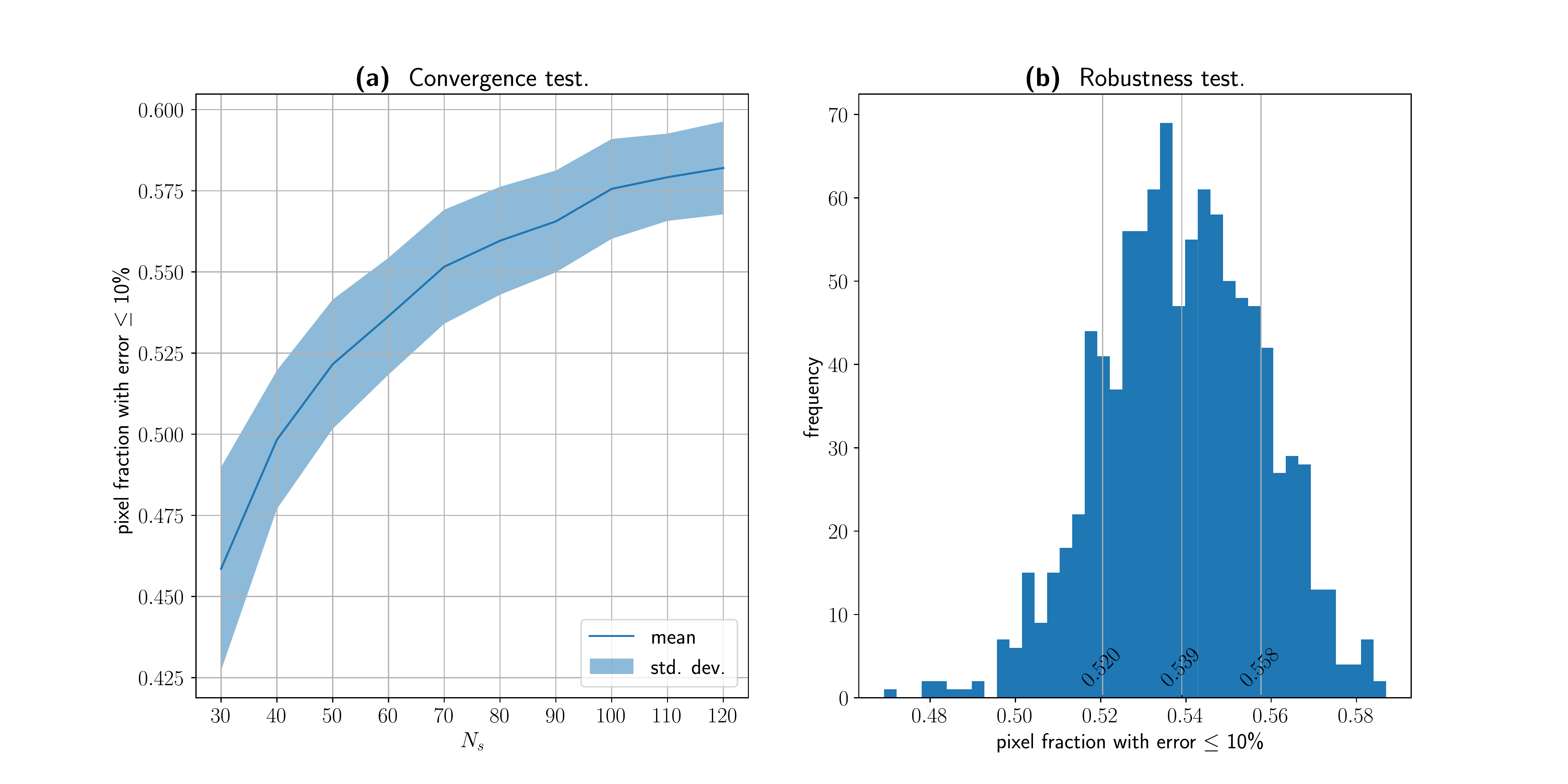}
\caption{\textbf{(a)} 
With $k = N_s$, 
we vary $N_s$ from 30 to 120, with a spacing of 10 units.  
For each value of $N_s$ we select 100 random uniform sets 
of views, evaluate the pixel fraction at 10\%
and report the mean and the standard deviation 
for each distribution.
The method 
delivers its best performance between 60 and 80 samples. 
In figure \textbf{(b)} a closer view of the pixel fraction 
distribution with respect the sample choice. 
With $N_s = k = 60$ we
pick 1000 random uniform sets of views, evaluate the 10\% 
pixel fraction for each set and plot the distribution. We
observe that the mean is 0.545 and the standard deviation is 
0.017, both the measures are expressed in terms of pixel 
fraction.}
\label{fig:c04-convergence-robustness}
\end{figure}

The convergence of the PA-POD method is studied in further detail in 
figure \ref{fig:c04-convergence-robustness}\textbf{(a)}. 
With $k = N_s$, 
we vary $N_s$ from 30 to 120, with a 
spacing of 10 units.  
For each value of $N_s$ we select 100 random uniform sets 
of views, evaluate the pixel fraction at 10\%
and report the mean and the standard deviation 
for each distribution.
PA-POD expresses its best performance where 
the sample is very 
sparse, in the range 30-60 samples. 
From 80 samples onwards the 
convergence reaches a plateau, this behavior is 
due to the non-interpolatory nature of the coefficients, 
see Discussion section for 
details. Furthermore we notice that the standard deviation is smaller than $0.02$ 
expressed in pixel fraction at 10\%, this implies that the 
performance is consistent with respect the views choice. 
In figure 
\ref{fig:c04-convergence-robustness}\textbf{(b)} we present a closer view of the pixel fraction distribution with respect the views choice. 
We consider the pixel fraction at 10\% error and we evaluate it for 
1000 uniformly random sets of sixty views, with $N_s = k=60$. 
The mean is 0.539 and the standard deviation is 
0.019 expressed in pixel fraction. This confirms the method 
is robust with respect the views selection. 
In both figures \ref{fig:c04-convergence-robustness}\textbf{(a)} and \textbf{(b)} we see that the mean values for the pixel fraction at the 10\% error are slightly smaller compared to the values in figure \ref{fig:c04-pod-relative-error-report}, 0.539, compared to 0.57. Evaluating the error pixel by pixel, repeating and collecting the median, as in figure \ref{fig:c04-pod-relative-error-report}, has a minimal regularization effect on the error metric, without any consequence for our overall conclusions. 

\begin{figure}[h!]
\centering
\includegraphics[width=\linewidth]{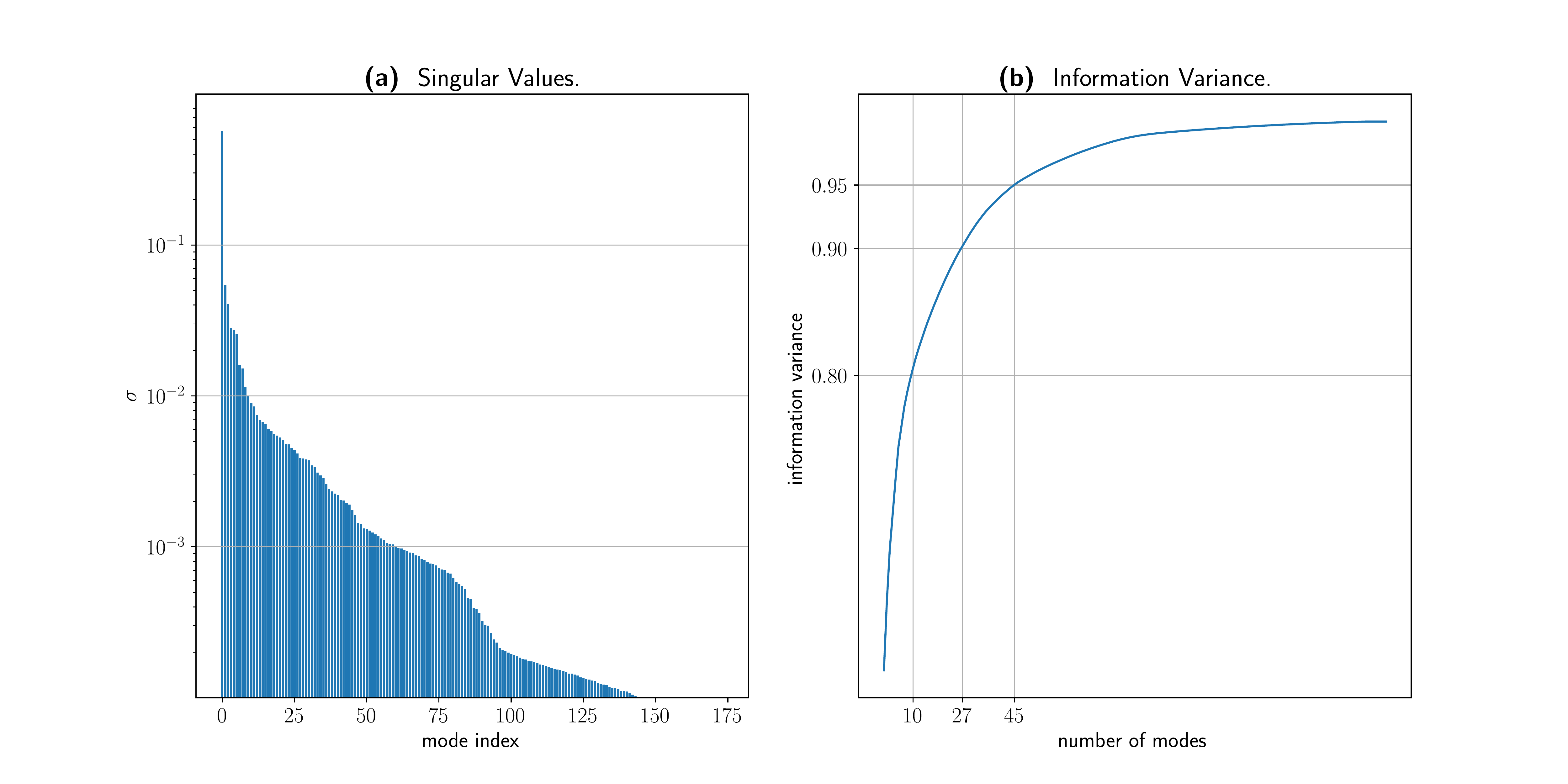}
\caption{
\textbf{(a)}
Logarithmic plot of the singular values spectrum for the IAEA PWR fuel 
assembly sinogram, and 
\textbf{(b)} the related information variance. 
We highlight that 
it takes respectively 8, 27, and 45 modes 
to capture 80\%, 90\% and 95\%
of the total information variance. 
}
\label{fig:c04-energy-spectrum}
\end{figure}

Figure \ref{fig:c04-energy-spectrum} 
plots the singular values for the 
PWR ground truth sinogram \textbf{(a)} and the corresponding information 
variance \textbf{(b)}. Figure \ref{fig:c04-energy-spectrum}\textbf{(b)}
highlights that a reduced set of modes 
captures the majority of the information variance. 
In particular 8, 27, and 45 modes 
capture 80\%, 90\% and 95\% 
of the total information variance respectively.
This result confirms the underlying assumption of our work: 
a singoram, experimentally measured or numerically computed,
is described by a limited set of dominant patterns,
that carefully combined together can provide a 
complete solution 
with a sparse set of angular views. 

\section*{Discussion}
\label{sec:pod}

Our approach was in part inspired by a particular flavor of Proper Orthogonal Decomposition that uses interpolation to estimate the coefficients, namely POD with Interpolation (PODI).  In this domain of interest, the research effort is devoted to finding the most accurate interpolator for each applicable case \cite{TezzeleDemoRozza2019MARINE, OrtaliDemoRozza2020MINE, DemoTezzeleRozza2019CRM, certified_rozza_2015}. 
PODI is particularly appealing because it is agnostic of the underlying physics of the phenomenon, it has no notion of the mathematical method describing the prototype, in ROM literature these methods are called non-intrusive \cite{certified_rozza_2015}. In this section we will briefly summarize the steps that define the method, we will show its limitations in the PGET case, and how we addressed them while developing PA-POD. 

In PODI the interpolatory coefficients are estimated by projecting the database matrix onto the POD space: $\mathbf{C} = \mathbf{U^*}\, \mathbf{\hat{S}}$. The  columns of $\mathbf{C}_{k\times N_s}$ construct a set of pairs $\{\mathbf{c}_j,\theta_j\}$, with $j = 0, \ldots, N_s-1$. Consider $\tilde\theta$ a rotation angle not included in our sample, we interpolate the coefficients $c_{ij}$ to get $\{\tilde{\mathbf{c}},\tilde{\theta}\}$ and reconstruct: 

\[
\mathbf{\tilde{s}}(y,\tilde{\theta})= \sum_{i=0}^{k-1} \mathbf{u}_i(y) \cdot \tilde{\mathbf{c}}. 
\]

$\mathbf{\tilde{C}}_{k\times 360} = \mathcal{I}(\mathbf{C}_{k \times N_s})$ is the coefficients matrix for the full 360 degrees, it is the result of the interpolation operator $\mathcal{I}$ applied to $\mathbf{\hat{S}}$ projected into the POD space.

\begin{figure}[h!]
\centering
\begin{tikzpicture}[x=0.75pt,y=0.75pt,yscale=-1,xscale=1]

\draw (300,50) node [anchor=north ][inner sep=0.75pt] [align=center] {
\large{\textbf{Physics-Aware POD}}};

\draw  [draw opacity=0][fill={rgb, 255:red, 255; green, 127; blue, 14 }  ,fill opacity=.3 ] (180,74) -- (420,74) -- (420,132) -- (180,132) -- cycle ;

\draw (300,80) node [anchor=north ][inner sep=0.75pt]  [align=center] {
\textbf{Define}\\
Real-Time Approximated Forward Model,\\
as the operator $\mathcal{R}(\theta)$.};

\draw  [draw opacity=0][fill={rgb, 255:red, 188; green, 189; blue, 34 }  ,fill opacity=.3 ] (200,138) -- (400,138) -- (400,205) -- (200,205) -- cycle ;
\draw (300,150) node [anchor=north ][inner sep=0.75pt]  [align=center]  {
\textbf{Evaluate} the real time model\\ for all the parameters:\\
$\mathbf{R}_{N\times 360} = \mathcal{R}([ 0^{\circ } ,\dotsc ,\ 359^{\circ }])$};


\draw (50,50) node [anchor=north][inner sep=0.75pt]   [align=center] {\large{\textbf{POD with Interpolation}}};

\draw  [draw opacity=0][fill={rgb, 255:red, 45; green, 160; blue, 44}  ,fill opacity=.3 ] (-50,80) -- (150,80) -- (150,158) -- (-50,158) -- cycle ;
\draw (50,90) node [anchor=north][inner sep=0.75pt]   [align=center] {
\textbf{Project} the database matrix\\
in the POD space,\\
$\mathbf{C}_{k\times N_s} = \mathbf{U^*}_{k\times N}\, \mathbf{S}_{N\times N_s},$\\
to construct $N_s$ pairs $\{\mathbf{c}_i,\theta_i\}_{N_s}$.
};

\draw  [draw opacity=0][fill={rgb, 255:red, 255; green, 127; blue, 14 }  ,fill opacity=.3 ] (-33,174) -- (133,174) -- (133,198) -- (-33,198) -- cycle ;
\draw (50,180) node [anchor=north][inner sep=0.75pt]   [align=center] {\textbf{Define} an interpolator: $\mathcal{I}$};

\draw  [draw opacity=0][fill={rgb, 255:red, 188; green, 189; blue, 34 }  ,fill opacity=.3 ] (-50,212) -- (150,212) -- (150,292) -- (-50,292) -- cycle ;
\draw (50,220) node [anchor=north][inner sep=0.75pt]  [align=center]  {
\textbf{Evaluate}, via the interpolator,\\
 the coefficients for all\\ the parameters:\\
$\tilde{\mathbf{C}}_{k\times 360} =\mathcal{I}\left(\mathbf{C}_{k\times N_s}\right)$};

\draw  [draw opacity=0][fill={rgb, 255:red, 45; green, 160; blue, 44}  ,fill opacity=.3 ] (200,212) -- (400,212) -- (400,292) -- (200,292) -- cycle ;

\draw (300,220) node [anchor=north][inner sep=0.75pt]    [align=center]{
\textbf{Project} the Real-Time Model\\ into the POD space\\
to obtain the coefficients matrix:\\
$\tilde{\mathbf{C}}_{k\times 360} =\mathbf{U^{*}{}}_{k\times N}\mathbf{R}_{N\times 360}$};


\draw  [draw opacity=0][fill={rgb, 255:red, 32; green, 119; blue, 180  }  ,fill opacity=.3 ] (50,312) -- (300,312) -- (300,375) -- (50,375) -- cycle ;

\draw (175,320) node [anchor=north][inner sep=0.75pt] [align=center]   {
The approximation is a combination\\
of coefficients and basis functions\\
$\tilde{\mathbf{S}}_{N\times 360} =\mathbf{U}_{N\times k}\tilde{\mathbf{C}}_{k\times 360}$};

\end{tikzpicture}
\caption{This picture shows a side by side comparison between the PODI and 
PA-POD. In PODI we interpolate accurate but scarce data, 
while with PA-POD we evaluate dense real time approximated data.}
\label{fig:papod_poi}
\end{figure}
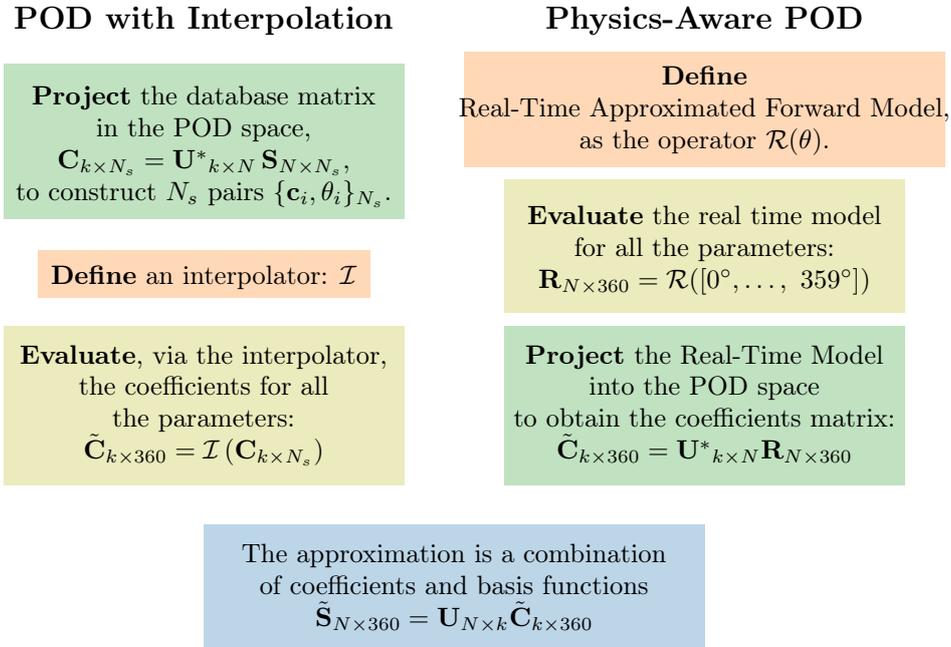

We compare the PODI and PA-POD approaches side by side in figure \ref{fig:papod_poi}. In both cases we look for the same solution structure:

\[
\mathbf{\tilde{S}}_{N\times 360} =\mathbf{U}_{N\times k}\tilde{\mathbf{C}}_{k\times 360}
\]

In both cases, the coefficients are evaluated by projecting available data into the POD space:

\[
\mathbf{C}_{k\times N_s}=\mathbf{U^*}_{k\times N}\,\mathbf{\hat{S}}_{N\times N_s} \quad \mathrm{vs.} \quad \tilde{\mathbf{C}}_{k\times 360} =\mathbf{U^{*}{}}_{k\times N}\mathbf{R}_{N\times 360},
\]
while PODI interpolates accurate but scarce 
data $\mathbf{\hat{S}}_{N\times N_s}$, 
conversely, PA-POD relies on approximate but 
dense data provided by the Real-Time Model $\mathbf{R}_{N\times 360}$.

\begin{figure}[!ht]
\centering
\includegraphics[width=\linewidth]{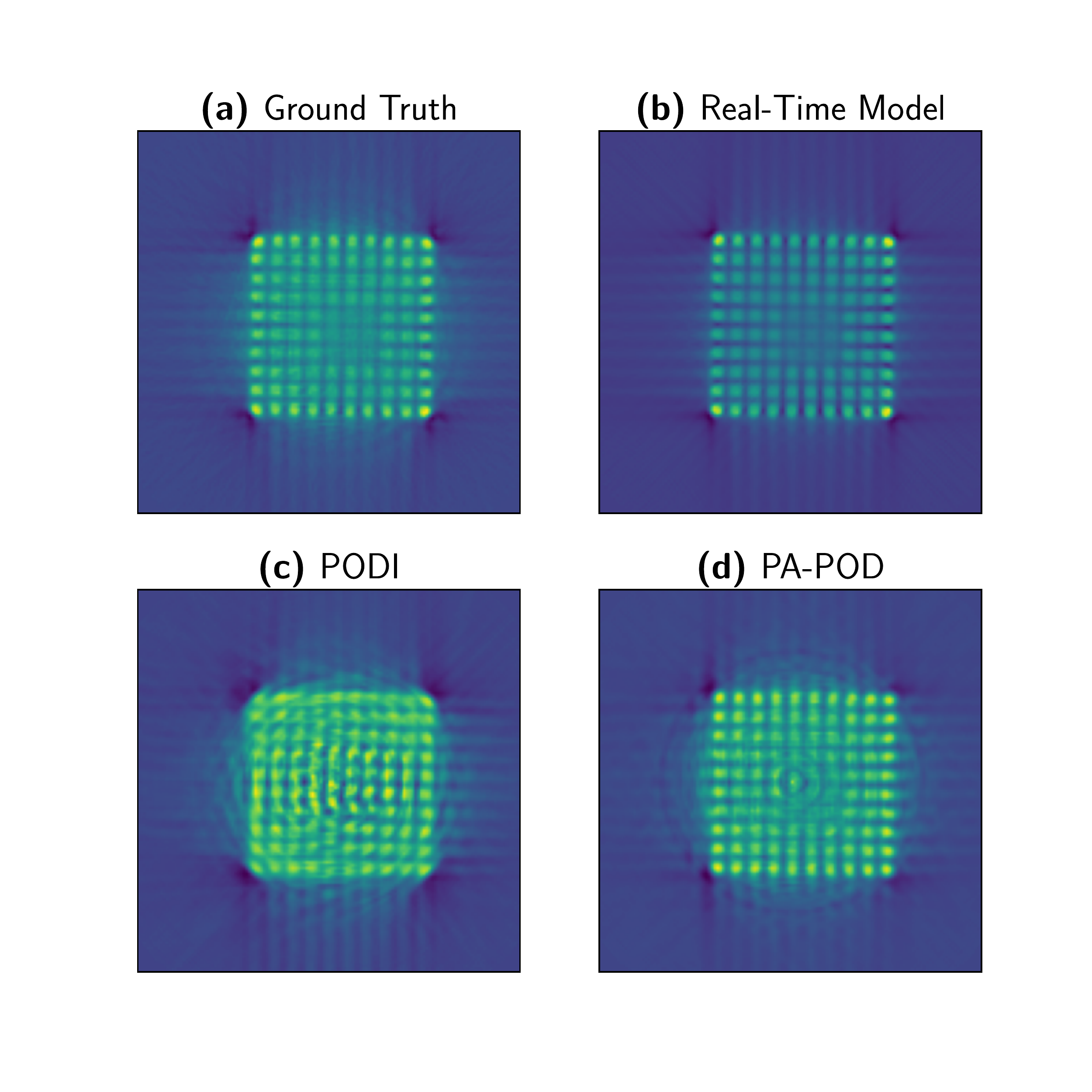}
\caption{Filtered backprojection for the IAEA data \textbf{(a)}, and the approximations we are studying, $N_s = k = 60$. 
Real-Time Model \textbf{(b)} and PA-POD \textbf{(d)} preserve the geometrical structure of the spent fuel, while PODI cannot provide a representative reconstruction \textbf{(c)}.}
\label{fig:c04-pod-example-iradon}
\end{figure}

In figure \ref{fig:c04-pod-example-iradon} we compare the FBP for several approximation techniques with the ground truth, \textbf{(a)}. 
The Real-Time Model \textbf{(b)} is an effective approximation in describing the assembly geometry. Figures \textbf{(c)} and \textbf{(d)} are constructed with $N_s = k = 60$. PODI fails to capture the structure of the fuel assembly \textbf{(c)}, while PA-POD preserves its geometrical structure \textbf{(d)} and recovers part of the secondary effects in the ground truth. 
Ring artifact is notable 
in the PA-POD reconstruction,
figure \ref{fig:c04-pod-example-iradon}\textbf{(d)}, 
this effect can be removed using the algorithm  
provided in reference \cite{Vo:18}.

\begin{figure}[h!]
\centering
\includegraphics[width=\linewidth]{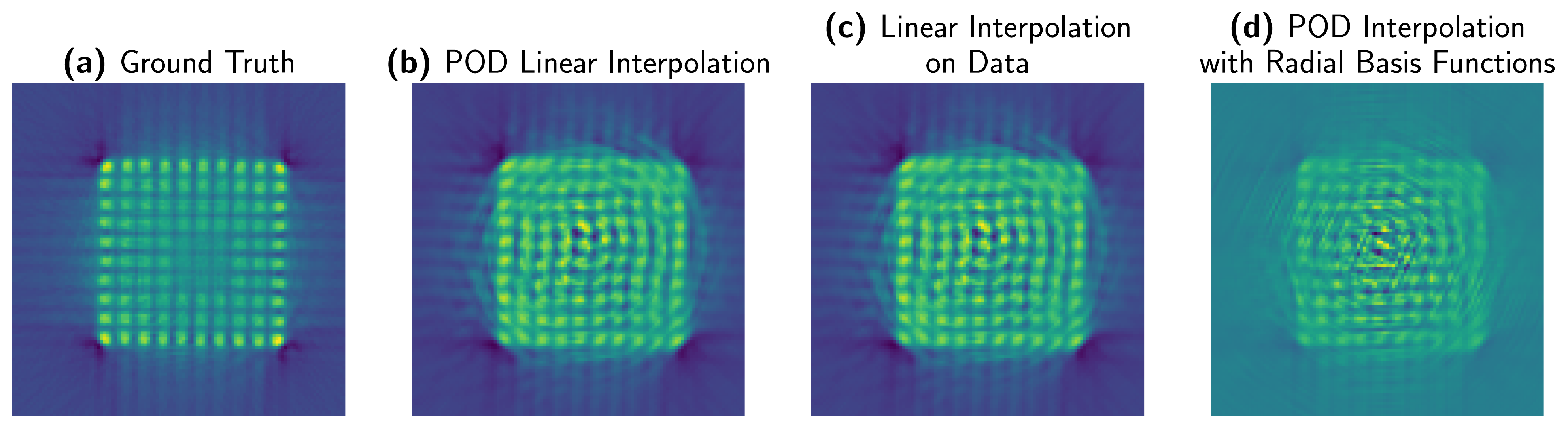}
\caption{This picture presents the filtered backprojection applied to three different interpolations, with $N_s =60$ and $k =30$. 
\textbf{(a)} Ground truth; 
\textbf{(b)} POD with linear interpolation of the coefficients; \textbf{(c)} Linear interpolation on the data; \textbf{(d)} POD with radial basis functions interpolation. The geometrical structure of the fuel assembly is severely compromised in all of three approximations. None of them preserves the assembly pin structure.}
\label{fig:interp-iradon}
\end{figure}

We tested three different interpolators. The results are reported in figure 
\ref{fig:interp-iradon}, the ground truth in figure \textbf{(a)}, 
\textbf{(b)} linear interpolation on the coefficients in the POD space, 
\textbf{(c)} linear interpolation on the data, 
\textbf{(d)} radial basis functions interpolation on the coefficients, 
but none could preserve the assembly geometry for $N_s = k = 60$.
In this phase of the exploration we relied on the EZyRB package 
developed at Scuola Internazionale Superiore Studi Avanzati \cite{demo18ezyrb}.

\begin{figure}[h!]
\centering
\includegraphics[width=\linewidth]{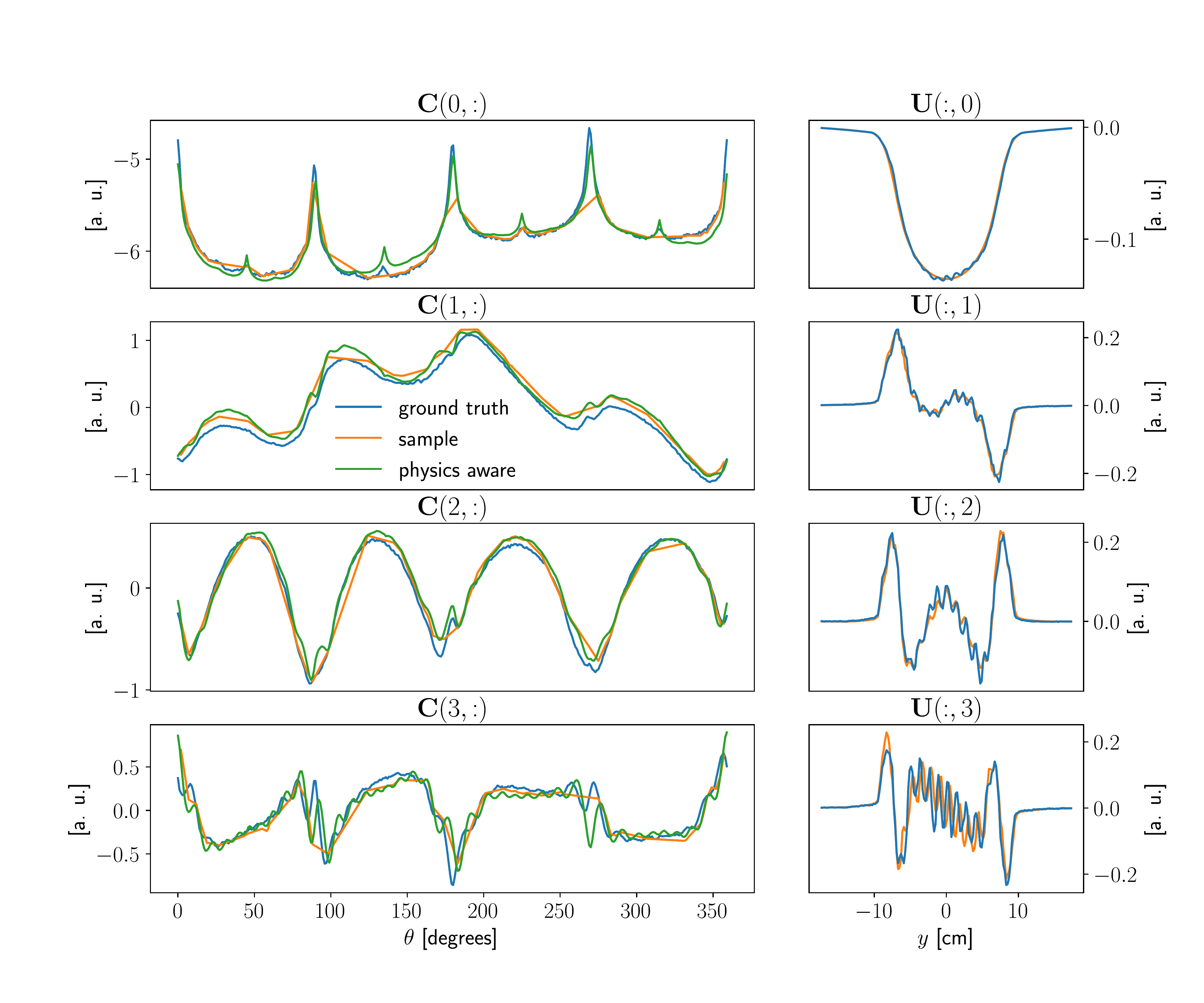}
\caption{Plots of the first four coefficients $\mathbf{C}(i,:)$ and modes $\mathbf{U}(:,i)$ for the PWR sinogram, with $k = N_s = 60$. A few 
discrepancies are visible in the sampled coefficients, 
while the physics-aware ones have a better definition especially in areas
characterized by high oscillations. 
The sampled modes, obtained applying SVD to $\mathbf{\hat{S}}$, 
are in good agreement with ones obtained using the whole dataset 
$\mathbf{S}$.}
\label{fig:c04-coeffs-overall}
\end{figure}

\begin{figure}[h!]
\centering
\includegraphics[width=\linewidth]{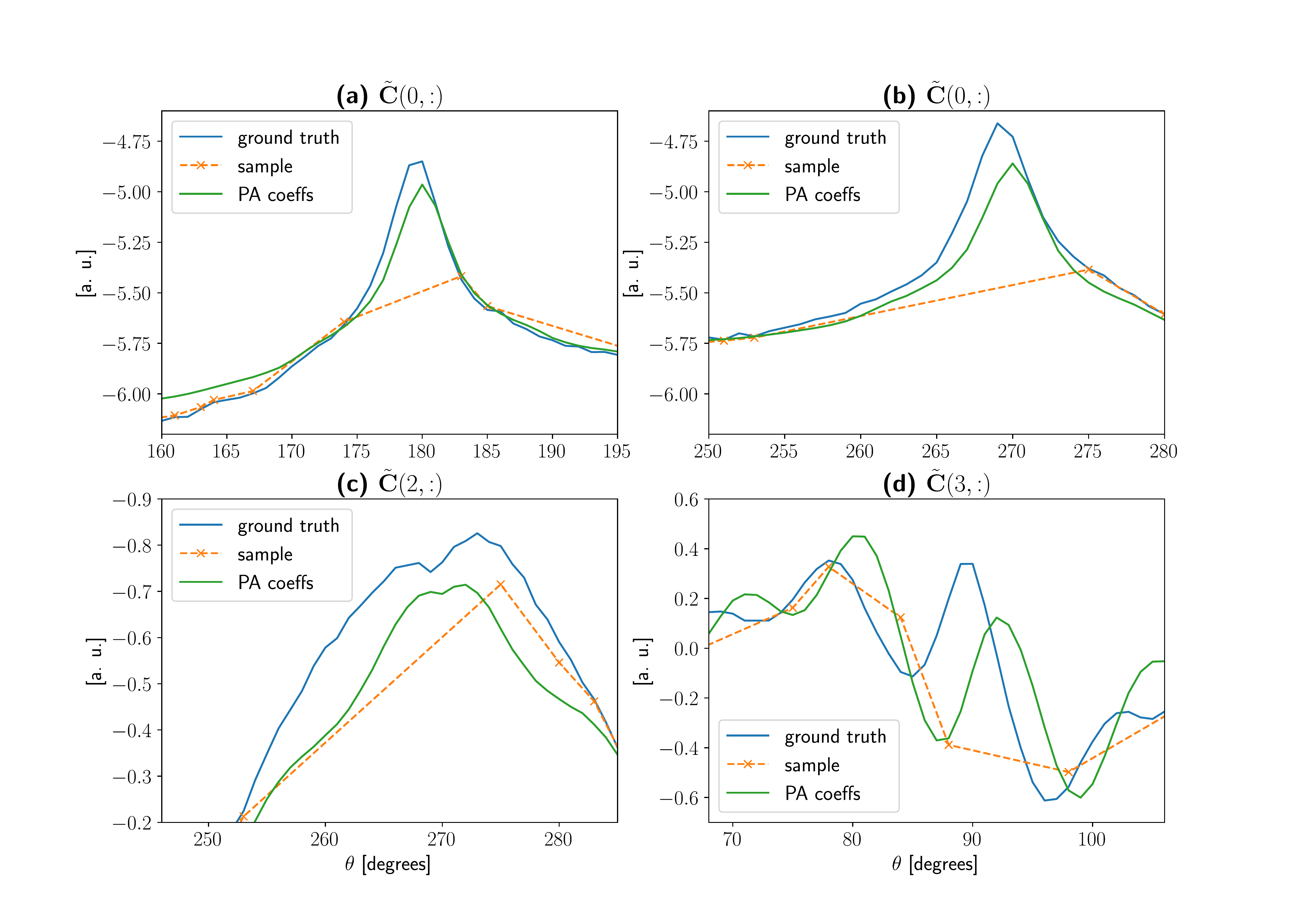}
\caption{In these plots a close up view of the ground truth coefficients, 
the sampled ones, and the physics-aware ones are reported. 
Given $N_s=k=60$  with randomly sampled 
views, we observed that the sampled coefficients do not capture enough 
data to account for the local structure. The information dropped during the 
sampling is in fact irrecoverable, unless we provide a physics-aware way to 
reconstruct the coefficients. It is this notion of ``physics awareness'' 
coefficients that inspired the main idea in this 
article.}
\label{fig:c04-coeffs-detail}
\end{figure}

The most important reason why interpolating methods show an inadequate performance is 
the sampling scarcity. In particular, figures \ref{fig:c04-coeffs-overall} and 
\ref{fig:c04-coeffs-detail}, show that the interpolation points fall too far apart to be 
effectively interpolated. 

Figure \ref{fig:c04-coeffs-overall} gives an overview of modes and coefficients, with 
$k = N_s = 60$. 
It shows the first four coefficients evaluated in three different ways. 
First, we project the ground truth into the POD space:
\[
\mathbf{C}_{\mathrm{ground\ truth}} = \mathbf{U}^*\, \mathbf{S},
\]
this is the best possible approximation of the expected coefficients. 
Second, we project the sample $\mathbf{\hat{S}}$ into the POD space, 
as prescribed by the PODI approach. 
Third, we evaluate the coefficients with the PA-POD method, 
projecting into the POD space an inexpensive approximation 
of the full sinogram, namely the result of the Real-Time Model.

We notice that while the sample preserves an overall fitting with the ground truth, 
important details are missed, causing the failure of the geometrical reconstruction. 
On the other hand, PA-POD coefficients are not interpolatory as in the PODI case but 
provide a better overall description of the assembly geometry. 
We notice that the modes $\mathbf{U}(:,i)$ for the database matrix $\mathbf{\hat{S}}$ 
are in good agreement with the ones obtained from the full sinogram $\mathbf{S}$, 
figure \ref{fig:c04-coeffs-overall}.

Figure \ref{fig:c04-coeffs-detail} shows several details of the coefficient matrices. 
The sampling is not sufficiently detailed to capture the actual shape of the 
coefficient, regardless of the accuracy of the interpolator, there is no way to recover 
information that is not in the data unless physics awareness of the prototype is 
restored, as we did with PA-POD.

It is evident, as depicted in figure \ref{fig:c04-coeffs-detail}, 
that PA-POD coefficients are not interpolatory. 
As a consequence, the resulting sinogram is not exactly the original one 
at the sampled locations.

\subsection*{Future Directions}

The ability of our PA-POD to create extended simulations at a low computational cost enables the construction of a digital twin, where actual measurements and simulated data are combined with a deep learning algorithm in a continuous update as measurements and simulations are produced in real-time. The benefit of doing so is to provide an integrated system, hardware, and digital twin for real-time high-accuracy verification of irradiated nuclear fuel assemblies using PGET. Our proposed approach can be readily adapted and used for other applications in nuclear safeguards and beyond, such as the creation of a dataset for the Tomographic Gamma Scanner technique \cite{VENKATARAMAN2007375}, which is a non-destructive assay technique, as well as in other radiography/tomography approaches for nuclear material and global security applications \cite{doi:10.1063/5.0075960}.

\section*{Methods}
\label{sec:methods}

\subsection*{Reference data used in the paper (IAEA Competition)}

Our method’s development and the associated verification have benefited from high-fidelity data from real case scenarios. The dataset provided by the International Atomic Energy Agency on the occasion of IAEA Tomographic and Analysis Challenge fit the purpose. The data are publicly available \cite{iaeauniteweb, iaeafullreport, 9059735}. Among several mockup data, few measurements are provided that were used in the paper. In the paper, the original names of the IAEA data, ``competition 3'', ``competition 4'', and ``competition 5'' are renamed into VVER, PWR, BWR respectively. In our presented results, we focused on the PWR case, 600-700 keV energy-deposition window, which it contains the full energy peak from the \ce{^{137}Cs} gamma emission at 661.7 keV, a major gamma emission of the spent nuclear fuel \cite{FAVALLI2016102}.

\subsection*{Real-Time Approximate Forward Model}

\begin{figure}[ht!]
\centering
\input{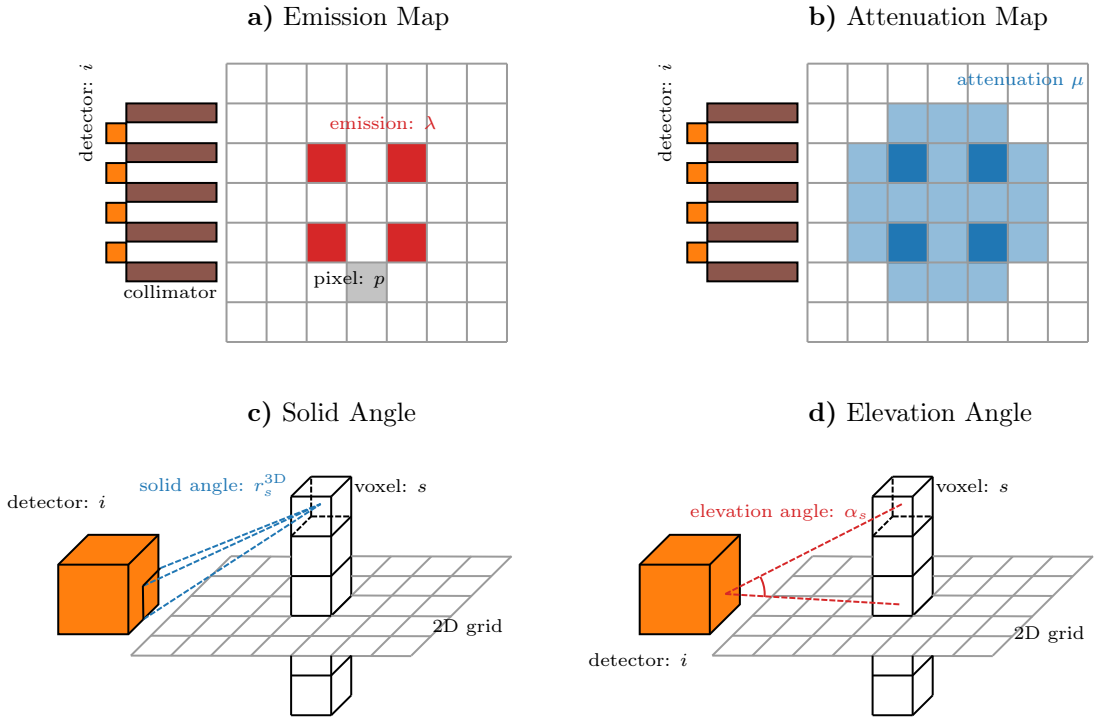}
\caption{A schematic view of the major players in the real time model implementation.}
\label{fig:pget-finland-sketch}
\end{figure}

Our Real-Time Model is built on the model of Backholm et al. \cite{Backholm2020}. A schematic representation of the model set is depicted in figure \ref{fig:pget-finland-sketch}. The investigation domain $\Omega \subset{\mathbb{R}^2}$ is the axial cross-section of a nuclear fuel assembly and we discretize it as a two-dimensional mesh. At each pixel $p$ of the discretized domain, we define $\lambda\ \in \mathbb{R}^{N_{\mathrm{pix}}}$ the emission values and the attenuation values $\mu\ \in \mathbb{R}^{N_{\mathrm{pix}}}$, $N_{\mathrm{pix}}$ is the number of pixels. We consider the values of $\lambda$ and $\mu$ constant over each pixel. At each detector $i$ we model the flux using a Lambert-Beer (exponential) attenuation model \cite{Backholm2020}:

\[
F(\lambda,\mu)_i = H(\mu)_{i,p}\, \lambda_{p} = \left[ r_{i,p}\,\mathrm{exp}\left(-c_{i,p}\, d^T_{i,p}\, \mu\right)\right]\, \lambda_{p}.
\]

The 360 views are obtained by rotating the pin positions, and so the assembly, with respect to the detectors. The emission and attenuation maps are consequently updated. The two-dimensional grid is further discretized, dividing the half-space above and below each pixel, in voxels.

Given this setting, the following quantities are defined. 

\begin{itemize}
    \item $r_{i,p}$ is the response of the detector $i$ with respect to the pixel $p$. Given a three-dimensional point, the barycenter of a voxel $s$, $r^{3D}_s$ is the solid angle between the voxel barycenter and the detector face, divided by $4\pi$. Follows $r_{i,p}$ is the average of the responses for all the voxels, insisting on the same pixel: $r_{i,p} = \frac{1}{N_{p,\mathrm{vox}}}\sum_{s=1}^{N_{p,\mathrm{vox}}} r^{3D}_s$.

    \item $d^T_{i,p}\, \mu$ is the integral of the discrete attenuation values on the segment connecting the pixel $p$ and the detector $i$. In our implementation we explicitly find the intersection between the connecting segment and the grid, being the attenuation piecewise constant, a zero-order quadrature is sufficient to exactly evaluate the integral.

    \item $c_{i,p}$ is an attenuation correction factor, evaluated as: $c_{i,p} = \frac{1}{r_{i,p}}\sum_{s=1}^{N_{p,\mathrm{vox}}} \frac{r^{3D}_s}{\mathrm{cos}(\alpha_s)}$. $\alpha_s$ is the angle between the voxel barycenter, the detector and the two dimensional plane. It accounts for the response matrix $r_{i,p}$ and the correction term $c_{i,p}$ are purely geometrical measures and are not affected by the rotation of the fuel assembly. $d^T_{i,p}\, \mu$ instead depends on the attenuation map, which depends on attenuation distribution, that changes with the assembly rotation and has to be assembled for each rotation angle. 
    
\end{itemize}

We tested three different mesh sizes, while $\Delta z = 10\ \mathrm{mm}$ was the same in all the three cases, the other dimensions are $\Delta x = \Delta y = (2.5 \ \mathrm{mm}, 1.0 \ \mathrm{mm}, 0.5 \ \mathrm{mm})$. The mesh size does not affect the accuracy of our results, the coarsest mesh size can be chosen to minimize the computational cost. 

Values for emission $\lambda$ attenuation $\mu$ are defined in \cite{Backholm2020}. Emission is expressed in arbitrary units, and is 0 for water and 100 for the spent fuel. Attenuation is $0.1356\ \mathrm{mm}^{-1}$ for the spent fuel, and $0.0085\ \mathrm{mm}^{-1}$ for water. Virta et al.  in \cite{Virta2020} provide an extensive exploration of such parameters in experimental settings. An approximation of the Real-Time approximate forward model is that it does not model the downscattering of the gamma rays.

The computational cost for a full sinogram according to this ray-tracing scheme is a few minutes on a desktop machine, and the assembly of all the matrices is a type of ``embarrassingly parallel'' task \cite{10.5555/2505461}, it can be parallelized both on a ``per detector'' basis or on a ``per degree'' basis.

\newpage
\clearpage

\bibliography{sample}
\bibliographystyle{abbrv}

\end{document}